\documentclass[11pt]{scrartcl}
\usepackage{hyperref}
\usepackage{epsfig}
\usepackage{times}
\usepackage[numbers]{natbib}
\usepackage{amsmath}
\usepackage{amsfonts}
\usepackage[latin1]{inputenc}
\usepackage{url}

\providecommand{\email}[1]{E-mail: \href{mailto:#1}{\texttt{#1}}}
\providecommand{\dprod}{\! \cdot \!}%

\begin{document}
\pagestyle{myheadings} \markright{J.B. Almeida / Standard-model
symmetry\dots}
\title{Standard-model symmetry in complexified
spacetime algebra}
\author{Jos\'e B. Almeida\\
\emph{Universidade do Minho, Departamento de F\'isica,}\\
\emph{4710-057 Braga, Portugal.}\\
\email{bda@fisica.uminho.pt}}


\date{\today}

\maketitle
\begin{abstract}                
Complexified spacetime algebra is defined as the geometric
(Clifford) algebra of spacetime with complex coefficients,
isomorphic $\mathcal{G}_{1,4}$. By resorting to matrix
representation by means of Dirac-Pauli gamma matrices, the paper
demonstrates isomorphism between subgroups of CSTA and $SU(3)$. It
is shown that the symmetry group of those subgroups is indeed
$U(1) \otimes SU(2) \otimes SU(3)$ and that there are 4 distinct
copies of this group within CSTA.
\end{abstract}
\section{Introduction}
The application of geometric algebra $\mathcal{G}_{1,3}$, a.k.a.\
Clifford algebra $\mathcal{C}\ell_{1,3}$, to quantum mechanics was
initiated by David Hestenes,
\citet{Hestenes73,Hestenes75,Hestenes86}, who was also responsible
for the designation \emph{Spacetime algebra} (STA) for the
geometric algebra of spacetime; the Cambridge group has also
brought important contributions, \citet{Lasenby93, Doran93}.
However none of those authors tackled the problem of extending the
application to the Standard Model; this was done under a different
approach by another group in Canada, \citet{Trayling01}, using the
7-dimensional $\mathcal{G}_{4,3}$ mapped to $\mathcal{G}_7$
through the use of scalar imaginary $j$. The choice of
7-dimensional space to accommodate $SU(3)$ symmetry arises quite
naturally, \citet{Lounesto01}, although it should be clear that
this is an oversized dimensionality, since $SU(3)$ is a group of
$3 \times 3$ matrices, while 3-dimensional geometric algebra has a
$2 \times 2$ matrix representation and 4-dimensional geometric
algebra is represented with $4 \times 4$ matrices. One would then
expect the space dimensionality corresponding to $SU(3)$ to lie
somewhere between 3D and 4D, if that was a possibility.

A suitable basis for $SU(3)$ is provided by the $3 \times 3$
Gell-Mann matrices, \citet{Cottingham98}; here we use those
matrices multiplied by $j$;
\begin{alignat}{3}
    \label{eq:gellmann}
    \hat{\lambda}_1 =& \begin{pmatrix}
    0 & j & 0 \\ j & 0 & 0 \\ 0 & 0 & 0
    \end{pmatrix}, &~~
    \hat{\lambda}_2 =& \begin{pmatrix}
    0 & 1 & 0 \\ -1 & 0 & 0 \\ 0 & 0 & 0
    \end{pmatrix}, &~~
    \hat{\lambda}_3 =& \begin{pmatrix}
    j & 0 & 0 \\ 0 & -j & 0 \\ 0 & 0 & 0
    \end{pmatrix}, \nonumber \\
    \hat{\lambda}_4 =& \begin{pmatrix}
    0 & 0 & j \\ 0 & 0 & 0 \\ j & 0 & 0
    \end{pmatrix}, &~~
    \hat{\lambda}_5 =& \begin{pmatrix}
    0 & 0 & 1 \\ 0 & 0 & 0 \\ -1 & 0 & 0
    \end{pmatrix}, &~~
    \hat{\lambda}_6 =& \begin{pmatrix}
    0 & 0 & 0 \\ 0 & 0 & j \\ 0 & j & 0
    \end{pmatrix}, \nonumber \\
    \hat{\lambda}_7 =& \begin{pmatrix}
    0 & 0 & 0 \\ 0 & 0 & 1 \\ 0 & -1 & 0
    \end{pmatrix}, &~~
    \hat{\lambda}_8 =& \frac{1}{\sqrt{3}} \begin{pmatrix}
    j & 0 & 0 \\ 0 & j & 0 \\ 0 & 0 & -2j
    \end{pmatrix}; &~~&
\end{alignat}
a hat $(\hat{\ })$ is used to signify that we are dealing with
matrices and not with geometric algebra multivectors with
corresponding designations. The $\hat{\lambda}$ matrices satisfy
the commutation relations
\begin{equation}
    \label{eq:lambdacommute}
    \hat{\lambda}_a \hat{\lambda}_b- \hat{\lambda}_b
    \hat{\lambda}_a
     = -2 \sum_{c=1}^8
    f_{abc} \hat{\lambda}_c,
\end{equation}
where the $f_{abc}$ are the structure constants. The $f_{abc}$ are
odd in the interchange of any pair of indices, and the
non-vanishing elements are given by the permutations of $f_{123} =
1$, $f_{147} = f_{246} = f_{257} = f_{345} = f_{516} = f_{637} =
1/2$, $f_{458} = f_{678} = \sqrt{3}/2$.

Since STA and 4-dimensional geometric algebras in general can be
represented by $4 \times 4$ matrices, it is clear that, in terms
of dimensionality at least, these algebras are large enough to
contain $SU(3)$; we will show that this is true if complex
coefficients are allowed.

The group of gauge symmetries is generally described by $U(1)
\otimes SU(2) \otimes SU(3)$, \citet{Trayling01}. Just as $SU(3)$
is associated to the Gell-Mann matrices, $SU(2)$ can be
represented by the $2 \times 2$ Pauli matrices
\begin{equation}
    \label{eq:pauli}
    \hat{\sigma}_1 = \begin{pmatrix}
    0 & 1 \\ 1 & 0 \end{pmatrix},~~
    \hat{\sigma}_2 = \begin{pmatrix}
    0 & -j \\
    j & 0 \end{pmatrix},~~
    \hat{\sigma}_3 = \begin{pmatrix}
    1 & 0 \\
    0 & -1 \end{pmatrix},
\end{equation}
which are also associated with the representation of 3-dimensional
geometric algebra $\mathcal{G}_3$. $U(1)$ is just the algebra of
complex numbers. The complexified spacetime algebra (CSTA)
contains the gauge symmetry group, avoiding the need to use extra
dimensions, as we show below. CSTA is 5-dimensional and can be
mapped to other 5-dimensional algebras, although we prefer not to
do so in order to facilitate geometrical interpretation.
\section{Complexified spacetime algebra}
Following \citet{Doran93}, with the notation conventions of paper
\citet{Gull93}, STA is defined as the geometric algebra of
Minkowski spacetime and is generated by a frame of orthonormal
vectors $\{\gamma_\mu\}$, $\mu = 0 \ldots 3$, that satisfy the
Dirac algebra
\begin{equation}
    \label{eq:diralgbr}
    \gamma_\mu \dprod \gamma_\nu = \frac{1}{2}\left(\gamma_\mu
    \gamma_\nu + \gamma_\nu \gamma_\mu \right) =
    \mathrm{diag}\left(+ - - -\right),
\end{equation}
but are to be considered as four independent vectors instead of
matrices. The full STA is spanned by the basis
\begin{equation}
    \label{eq:stabasis}
    1,~~~~ \left\{\gamma_\mu \right\},~~~~ \left\{\sigma_k, i
    \sigma_k \right\},~~~~ \left\{i \gamma_\mu \right\},~~~~i.
\end{equation}
Here $i\equiv \gamma_0\gamma_1\gamma_2\gamma_3$ is the unit
pseudoscalar; it anticommutes with vectors and trivectors and
$i^2=-1$. The spacetime bivectors $\left\{\sigma_k\right\}$, $k=1
\ldots 3$ are defined by
\begin{equation}
    \label{eq:sigmaequiv}
    \sigma_k \equiv \gamma_k \gamma_0,
\end{equation}
and represent an orthonormal frame of 3-dimensional space. The
$\{\sigma_k\}$ generate the Pauli algebra of space, so that
relative vectors $a^k \sigma_k$ are viewed as spacetime bivectors.
We will assume that the author is familiar with STA operations and
will not go into details in this work; consultation of the cited
references should clarify any possible doubts. CSTA doubles STA by
intervention of the scalar imaginary so that the coefficients
multiplying the basis elements in a multivector can be complex
numbers. CSTA is the complex $\mathcal{G}_{1,3}$ algebra,
isomorphic to the real 5-dimensional algebra $\mathcal{G}_{1,4}$,
but in this work we always use the complex 4-dimensional
alternative.

We follow \citet{Doran93} for the matrix representation of STA and
CSTA defining the $\hat{\gamma}$-matrices in the standard
Dirac-Pauli representation, \citet{Halzen84},
\begin{equation}
    \label{eq:diracpauli}
    \hat{\gamma}_0 = \begin{pmatrix}
      I & 0 \\
      0 & -I \
    \end{pmatrix},~~~~ \hat{\gamma}_k =  \begin{pmatrix}
      0 & -\hat{\sigma}_k \\
      ~\hat{\sigma}_k & ~0 \
    \end{pmatrix}.
\end{equation}
When needed we will also make use of $4 \times 4$
$\hat{\sigma}_\mu$ matrices defined by
\begin{equation}
    \label{eq:4sigma}
    \hat{\sigma}_0 = \hat{\gamma}_0,~~~~ \hat{\sigma}_k =
    \hat{\gamma}_k \hat{\gamma}_0,
\end{equation}
which can be seen as the matrix equivalent of relation
(\ref{eq:sigmaequiv}). We denote these $4 \times 4$ matrices in
the same way as their $2 \times 2$ counterparts; this should not
cause any confusion since the context will always make our intents
clear.

In order to discuss $SU(3)$ symmetry in CSTA we note that $4
\times 4$ counterparts to Gell-Mann's matrices (\ref{eq:gellmann})
can be built by inserting an extra column and an extra row of
zeroes at position 3; these matrices take the same designation as
the original ones, since we should always know the dimension
pertinent to each situation and confusion should never arise. The
eight $4 \times 4$ matrices are
\begin{alignat}{3}
    \label{eq:gellmann4}
    \hat{\lambda}_1 =& \begin{pmatrix}
    0 & j & 0 & 0 \\ j & 0 & 0 & 0 \\ 0 & 0 & 0 & 0 \\0 & 0 & 0 & 0
    \end{pmatrix}, &~~
    \hat{\lambda}_2 =& \begin{pmatrix}
    0 & 1 & 0 & 0 \\ -1 & 0 & 0 & 0 \\ 0 & 0 & 0 & 0 \\ 0 & 0 & 0 & 0
    \end{pmatrix}, &~~
    \hat{\lambda}_3 =& \begin{pmatrix}
    j & 0 & 0 & 0 \\ 0 & -j & 0 & 0 \\ 0 & 0 & 0 & 0 \\ 0 & 0 & 0 & 0
    \end{pmatrix}, \nonumber \\
    \hat{\lambda}_4 =& \begin{pmatrix}
    0 & 0 & 0 & j \\ 0 & 0 & 0 & 0 \\ 0 & 0 & 0 & 0 \\ j & 0 & 0 & 0
    \end{pmatrix}, &~~
    \hat{\lambda}_5 =& \begin{pmatrix}
    0 & 0 & 0 & 1 \\ 0 & 0 & 0 & 0 \\ 0 & 0 & 0 & 0 \\ -1 & 0 & 0 & 0
    \end{pmatrix}, &~~
    \hat{\lambda}_6 =& \begin{pmatrix}
    0 & 0 & 0 & 0 \\ 0 & 0 & 0 & j \\ 0 & 0 & 0 & 0 \\ 0 & j & 0 & 0
    \end{pmatrix}, \nonumber \\
    \hat{\lambda}_7 =& \begin{pmatrix}
    0 & 0 & 0 & 0 \\ 0 & 0 & 0 & 1 \\ 0 & 0 & 0 & 0 \\ 0 & -1 & 0 & 0
    \end{pmatrix}, &~~
    \hat{\lambda}_8 =& \frac{1}{\sqrt{3}} \begin{pmatrix}
    j & 0 & 0 & 0 \\ 0 & j & 0 & 0 \\ 0 & 0 & 0 & 0 \\ 0 & 0 & 0 &
    -2j
    \end{pmatrix}. &~~&
\end{alignat}
These matrices are obviously a basis for $SU(3)$ in the algebra of
4-dimensional matrices, in view of the way they were built. We can
now try to find the geometric equations defining the equivalents
of the $\hat{\lambda}$ matrices in CSTA; we do this by defing the
four basic elements
\begin{equation}
\begin{split}
    \label{eq:su3csta1}
    \lambda_1 = i (\gamma_1 - \sigma_1) / 2, &~~
    \lambda_2 = i (\gamma_2 - \sigma_2) / 2, \\
    \lambda_4 = (j \sigma_1 + \gamma_2) / 2, &~~
    \lambda_5 = (j \sigma_2 - \gamma_1) / 2;
\end{split}
\end{equation}
the remaining 4 elements can be defined by recursion
\begin{equation}
\begin{split}
    \label{eq:su3csta2}
    \lambda_3 = - \lambda_1 \lambda_2, &~~
    \lambda_6 = \lambda_1 \lambda_5 - \lambda_5 \lambda_1 \\
    \lambda_7 = \lambda_5 \lambda_2 - \lambda_2 \lambda_5, &~~
    \lambda_8 = (-2 \lambda_4 \lambda_5 - \lambda_3)/\sqrt{3}.
\end{split}
\end{equation}

The $\lambda$-multivectors defined by Eqs.\ (\ref{eq:su3csta1})
and (\ref{eq:su3csta2}) observe the same commutating relations as
their matrix equivalents; it is therefore clear that the CSTA
subalgebra provided by the $\lambda$-multivectors is isomorphic to
$SU(3)$.

Having demonstrated $SU(3)$ symmetry we need to demonstrate the
simultaneous existence of $U(1)$ and $SU(2)$. For the former it is
sufficient to note that the dual of Eqs.\ (\ref{eq:su3csta1}) and
(\ref{eq:su3csta2}) is obtained by simply multiplying every
$\gamma$-vector by the scalar imaginary $j$. Indeed, by
multiplying all $\gamma$-vectors by the same complex number we
preserve $SU(3)$ and thus prove the product $U(1) \otimes SU(3)$.
Then we note that the definitions allow permutation of the indices
$1, 2, 3$ in the $\gamma$-vectors, providing the desired $SU(2)$
symmetry. In appendix \ref{a1} we show the matrix equivalents
corresponding to these permutations.

Eqs.\ (\ref{eq:su3csta1}) define two pairs of anti-commuting
multivectors, $\lambda_1$--$\lambda_2$ and
$\lambda_4$--$\lambda_5$. The first pair produces the elements
$(1,2)$ and $(2,1)$ of the matrix representation, while the second
pair produces elements $(1,4)$ and $(4,1)$. Redefining the first
pair as
\begin{equation}
    \label{eq:1stpair}
    \lambda_1 = -i(\gamma_1 + \sigma_1)/2,~~ \lambda_2 = -i
    (\gamma_2 + \sigma_2)/2,
\end{equation}
we get elements $(3,4)$ and $(4,3)$, which are compatible with the
remaining equations, creating a copy of the gauge group. In a
similar fashion we can redefine the second pair as
\begin{equation}
    \label{eq:2ndpair}
    \lambda_4 = (j \sigma_1 - \gamma_2)/2,~~ \lambda_5 = j
    (\sigma_2 + \gamma_1)/2,
\end{equation}
which produces elements $(2,3)$ and $(3,2)$ of the matrix
representation. Combining the two alternatives for both pairs of
multivectors we get a total of 4 copies for the gauge group,
present in CSTA; in matrix representation the different copies are
characterized by the row/column that gets filled with zeros in the
0th permutation; the matrix representation for the remaining 3
copies is shown in appendix \ref{a2}. Each copy allows
multiplication by a complex number, yielding $U(1)$ symmetry, and
permutation of the $\gamma$-vector indices, resulting in $SU(2)$
symmetry.

Examination of the matrix representations in appendix \ref{a2}
shows that the 4 gauge group copies are not independent; in fact
only 4 matrices of copy $\#1$ are independent from those of copy
$\#0$, while $\hat{\lambda}_4$ and $\hat{\lambda}_5$ are common.
Proceeding to copy $\#2$, only two new independent matrices are
added, precisely $\hat{\lambda}_4$ and $\hat{\lambda}_5$, while
all matrices in copy $\#3$ can be obtained by linear combinations
from the previous ones; in total there are 16 independent matrices
as one could expect from the fact that $\mathcal{G}_{1,3}$ is a
graded 16-dimensional space. The set of 3 independent copies is
possibly connected to another $SU(2)$ symmetry.
\section{Conclusion}
The geometric algebra of spacetime (STA) is represented by $4
\times 4$ Dirac-Pauli matrices, while $SU(3)$ is the special group
of unitary $3 \times 3$ matrices. We introduced complex
coefficients in STA, defining the CSTA 5-dimensional algebra,
which shares the matrix representation of STA but allows each
matrix to be multiplied by a complex number. Within the matrix
representation of CSTA we were able to find an isomorphism of
$SU(3)$, which we could define in terms of multivectors and
geometric products.

$SU(3)$ isomorphism in CSTA is richer than the original $3 \times
3$ matrix one, first because the complex coefficients introduce
$U(1)$ symmetry and second because permutation among the frame
vectors present in the multivector definitions provide $SU(2)$
symmetry. The result is the standard-model $U(1) \otimes SU(2)
\otimes SU(3)$ symmetry group, expressed in Minkowski spacetime
provided with complex numbers, something familiar to most
physicists. We showed also that CSTA contains 4 inter-dependent
copies of the standard-model group from which 3 can be considered
independent of each other; this is possibly a manifestation
another $SU(2)$ symmetry.

We think it will be possible to build on the previous formalism to
reformulate the standard-model and gauge theories in terms of
geometric algebra, namely CSTA; some work is already under way in
that direction.
\appendix
\section{\label{a1}Matrix equivalents for $SU(2)$ permutations}
The matrices resulting from permutations of indices $1,2,3$ in the
$\gamma$-vectors can be quickly obtained by performing these
permutations in Eqs.\ (\ref{eq:su3csta1}) and then applying Eqs.\
(\ref{eq:su3csta2}) unchanged. The first right permutation of
$\gamma$-vector produces the equivalent matrices
\begin{alignat}{2}
    \label{eq:perm1}
    \hat{\lambda}_1 =& \begin{pmatrix}
    0 & 1 & 0 & 0 \\ -1 & 0 & 0 & 0 \\ 0 & 0 & 0 & 0 \\0 & 0 & 0 & 0
    \end{pmatrix}, &~~
    \hat{\lambda}_2 =& \begin{pmatrix}
    j & 0 & 0 & 0 \\ 0 & -j & 0 & 0 \\ 0 & 0 & 0 & 0 \\ 0 & 0 & 0 & 0
    \end{pmatrix}, \nonumber \\
    \hat{\lambda}_3 =& \begin{pmatrix}
    0 & j & 0 & 0 \\ j & 0 & 0 & 0 \\ 0 & 0 & 0 & 0 \\ 0 & 0 & 0 & 0
    \end{pmatrix}, \nonumber &~~
    \hat{\lambda}_4 =& \frac{1}{2}\begin{pmatrix}
    0 & 0 & -1& 1 \\ 0 & 0 & -1& 1\\
    1 & 1 & 0 & 0 \\ -1 & -1 & 0 & 0
    \end{pmatrix}, \\
    \hat{\lambda}_5 =& \frac{1}{2}\begin{pmatrix}
    0 & 0 & j & -j \\ 0 & 0 & j & -j \\
    j & j & 0 & 0 \\ -j & -j & 0 & 0
    \end{pmatrix}, &~~
    \hat{\lambda}_6 =& \frac{1}{2}\begin{pmatrix}
    0 & 0 & j & -j \\ 0 & 0 & -j & j \\
    j & -j & 0 & 0 \\ -j & j & 0 & 0
    \end{pmatrix}, \nonumber \\
    \hat{\lambda}_7 =& \frac{1}{2}\begin{pmatrix}
    0 & 0 & 1 & -1 \\ 0 & 0 & -1 & 1 \\
    -1 & 1 & 0 & 0 \\ 1 & -1 & 0 & 0
    \end{pmatrix}, &~~
    \hat{\lambda}_8 =& \frac{1}{\sqrt{3}} \begin{pmatrix}
    j & 0 & 0 & 0 \\ 0 & j & 0 & 0 \\ 0 & 0 & -j & j \\ 0 & 0 & j &
    -j
    \end{pmatrix}.
\end{alignat}
The second right permutation produces the matrices
\begin{alignat}{2}
    \label{eq:perm2}
    \hat{\lambda}_1 =& \begin{pmatrix}
    j & 0 & 0 & 0 \\ 0 & -j & 0 & 0 \\ 0 & 0 & 0 & 0 \\0 & 0 & 0 & 0
    \end{pmatrix}, &~~
    \hat{\lambda}_2 =& \begin{pmatrix}
    0 & j & 0 & 0 \\ j & 0 & 0 & 0 \\ 0 & 0 & 0 & 0 \\ 0 & 0 & 0 & 0
    \end{pmatrix}, \nonumber \\
    \hat{\lambda}_3 =& \begin{pmatrix}
    0 & 1 & 0 & 0 \\ -1 & 0 & 0 & 0 \\ 0 & 0 & 0 & 0 \\ 0 & 0 & 0 & 0
    \end{pmatrix}, \nonumber &~~
    \hat{\lambda}_4 =& \frac{1}{2}\begin{pmatrix}
    0 & 0 & j& -1 \\ 0 & 0 & -1& -j\\
    j & 1 & 0 & 0 \\ 1 & -j & 0 & 0
    \end{pmatrix}, \\
    \hat{\lambda}_5 =& \frac{1}{2}\begin{pmatrix}
    0 & 0 & 1 & j \\ 0 & 0 & j & -1 \\
    -1 & j & 0 & 0 \\ j & 1 & 0 & 0
    \end{pmatrix}, &~~
    \hat{\lambda}_6 =& \frac{1}{2}\begin{pmatrix}
    0 & 0 & j & -1 \\ 0 & 0 & 1 & j \\
    j & -1 & 0 & 0 \\ 1 & j & 0 & 0
    \end{pmatrix}, \nonumber \\
    \hat{\lambda}_7 =& \frac{1}{2}\begin{pmatrix}
    0 & 0 & 1 & j \\ 0 & 0 & -j & 1 \\
    -1 & -j & 0 & 0 \\ j & -1 & 0 & 0
    \end{pmatrix}, &~~
    \hat{\lambda}_8 =& \frac{1}{\sqrt{3}} \begin{pmatrix}
    j & 0 & 0 & 0 \\ 0 & j & 0 & 0 \\ 0 & 0 & -j & 1 \\ 0 & 0 & -1 &
    -j
    \end{pmatrix}.
\end{alignat}
\section{\label{a2}Multiple copies of the gauge group}
Equations (\ref{eq:su3csta1}) define one of the gauge group copies
coexisting in CSTA; by convention we label this as copy $\#0$;
copy $\#1$ is obtained by replacing the first of those equations
with Eq.\ (\ref{eq:1stpair}), resulting in the following matrices
\begin{alignat}{3}
    \label{eq:copy_1}
    \hat{\lambda}_1 =& \begin{pmatrix}
    0 & 0 & 0 & 0 \\ 0 & 0 & 0 & 0 \\ 0 & 0 & 0 & j \\0 & 0 & j & 0
    \end{pmatrix}, &~~
    \hat{\lambda}_2 =& \begin{pmatrix}
    0 & 0 & 0 & 0 \\ 0 & 0 & 0 & 0 \\ 0 & 0 & 0 & 1 \\ 0 & 0 & -1 & 0
    \end{pmatrix}, &~~
    \hat{\lambda}_3 =& \begin{pmatrix}
    0 & 0 & 0 & 0 \\ 0 & 0 & 0 & 0 \\ 0 & 0 & j & 0 \\ 0 & 0 & 0 &
    -j
    \end{pmatrix}, \nonumber \\
    \hat{\lambda}_4 =& \begin{pmatrix}
    0 & 0 & 0 & j \\ 0 & 0 & 0 & 0 \\ 0 & 0 & 0 & 0 \\ j & 0 & 0 & 0
    \end{pmatrix}, &~~
    \hat{\lambda}_5 =& \begin{pmatrix}
    0 & 0 & 0 & 1 \\ 0 & 0 & 0 & 0 \\ 0 & 0 & 0 & 0 \\ -1 & 0 & 0 & 0
    \end{pmatrix}, &~~
    \hat{\lambda}_6 =& \begin{pmatrix}
    0 & 0 & -j & 0 \\ 0 & 0 & 0 & 0 \\ -j & 0 & 0 & 0 \\ 0 & 0 & 0 & 0
    \end{pmatrix}, \nonumber \\
    \hat{\lambda}_7 =& \begin{pmatrix}
    0 & 0 & -1 & 0 \\ 0 & 0 & 0 & 0 \\ 1 & 0 & 0 & 0 \\ 0 & 0 & 0 & 0
    \end{pmatrix}, &~~
    \hat{\lambda}_8 =& \frac{1}{\sqrt{3}} \begin{pmatrix}
    2j & 0 & 0 & 0 \\ 0 & 0 & 0 & 0 \\ 0 & 0 & -j & 0 \\ 0 & 0 & 0 &
    -j
    \end{pmatrix}. &~~&
\end{alignat}
Similarly copy $\#2$ is obtained by replacing the second Eq.\
(\ref{eq:su3csta1}) with Eq.\ (\ref{eq:2ndpair})
\begin{alignat}{3}
    \label{eq:copy_2}
    \hat{\lambda}_1 =& \begin{pmatrix}
    0 & j & 0 & 0 \\ j & 0 & 0 & 0 \\ 0 & 0 & 0 & 0 \\0 & 0 & 0 & 0
    \end{pmatrix}, &~~
    \hat{\lambda}_2 =& \begin{pmatrix}
    0 & 1 & 0 & 0 \\ -1 & 0 & 0 & 0 \\ 0 & 0 & 0 & 0 \\ 0 & 0 & 0 & 0
    \end{pmatrix}, &~~
    \hat{\lambda}_3 =& \begin{pmatrix}
    j & 0 & 0 & 0 \\ 0 & -j & 0 & 0 \\ 0 & 0 & 0 & 0 \\ 0 & 0 & 0 & 0
    \end{pmatrix}, \nonumber \\
    \hat{\lambda}_4 =& \begin{pmatrix}
    0 & 0 & 0 & 0 \\ 0 & 0 & j & 0 \\ 0 & j & 0 & 0 \\ 0 & 0 & 0 & 0
    \end{pmatrix}, &~~
    \hat{\lambda}_5 =& \begin{pmatrix}
    0 & 0 & 0 & 0 \\ 0 & 0 & -1 & 0 \\ 0 & 1 & 0 & 0 \\ 0 & 0 & 0 & 0
    \end{pmatrix}, &~~
    \hat{\lambda}_6 =& \begin{pmatrix}
    0 & 0 & -j & 0 \\ 0 & 0 & 0 & 0 \\ -j & 0 & 0 & 0 \\ 0 & 0 & 0 & 0
    \end{pmatrix}, \nonumber \\
    \hat{\lambda}_7 =& \begin{pmatrix}
    0 & 0 & 1 & 0 \\ 0 & 0 & 0 & 0 \\ 1 & 0 & 0 & 0 \\ 0 & 0 & 0 & 0
    \end{pmatrix}, &~~
    \hat{\lambda}_8 =& \frac{1}{\sqrt{3}} \begin{pmatrix}
    -j & 0 & 0 & 0 \\ 0 & -j & 0 & 0 \\ 0 & 0 & 2j & 0 \\ 0 & 0 & 0 &
    0
    \end{pmatrix}. &~~&
\end{alignat}
Finally copy $\#3$ is obtained by making the two previous
replacements simultaneously
\begin{alignat}{3}
    \label{eq:copy_3}
    \hat{\lambda}_1 =& \begin{pmatrix}
    0 & 0 & 0 & 0 \\ 0 & 0 & 0 & 0 \\ 0 & 0 & 0 & j \\0 & 0 & j & 0
    \end{pmatrix}, &~~
    \hat{\lambda}_2 =& \begin{pmatrix}
    0 & 0 & 0 & 0 \\ 0 & 0 & 0 & 0 \\ 0 & 0 & 0 & 1 \\ 0 & 0 & -1 & 0
    \end{pmatrix}, &~~
    \hat{\lambda}_3 =& \begin{pmatrix}
    0 & 0 & 0 & 0 \\ 0 & 0 & 0 & 0 \\ 0 & 0 & j & 0 \\ 0 & 0 & 0 &
    -j
    \end{pmatrix}, \nonumber \\
    \hat{\lambda}_4 =& \begin{pmatrix}
    0 & 0 & 0 & 0 \\ 0 & 0 & j & 0 \\ 0 & j & 0 & 0 \\ 0 & 0 & 0 & 0
    \end{pmatrix}, &~~
    \hat{\lambda}_5 =& \begin{pmatrix}
    0 & 0 & 0 & 0 \\ 0 & 0 & -1 & 0 \\ 0 & 1 & 0 & 0 \\ 0 & 0 & 0 & 0
    \end{pmatrix}, &~~
    \hat{\lambda}_6 =& \begin{pmatrix}
    0 & 0 & 0 & 0 \\ 0 & 0 & 0 & j \\ 0 & 0 & 0 & 0 \\ 0 & j & 0 & 0
    \end{pmatrix}, \nonumber \\
    \hat{\lambda}_7 =& \begin{pmatrix}
    0 & 0 & 0 & 0 \\ 0 & 0 & 0 & -1 \\ 0 & 0 & 0 & 0 \\ 0 & 1 & 0 & 0
    \end{pmatrix}, &~~
    \hat{\lambda}_8 =& \frac{1}{\sqrt{3}} \begin{pmatrix}
    0 & 0 & 0 & 0 \\ 0 & -2j & 0 & 0 \\ 0 & 0 & j & 0 \\ 0 & 0 & 0 &
    j
    \end{pmatrix}. &~~&
\end{alignat}
  \bibliographystyle{unsrtbda}
  \bibliography{Abrev,aberrations,assistentes}   

\end{document}